\documentclass[12pt,leqno,fleqn]{amsart}
\usepackage{amsmath,amstext,amsthm,amssymb,amsxtra}
\usepackage[top=1.5in, bottom=1.5in, left=1.25in, right=1.25in]	{geometry}
\usepackage[normalem]{ulem}
\usepackage{txfonts} 
\usepackage[T1]{fontenc}
\usepackage{lmodern}
\usepackage{tikz}\usetikzlibrary{arrows}

\usepackage{euler}   

\makeatletter
\DeclareRobustCommand\widecheck[1]{{\mathpalette\@widecheck{#1}}}
\def\@widecheck#1#2{%
    \setbox\z@\hbox{\m@th$#1#2$}%
    \setbox\tw@\hbox{\m@th$#1%
       \widehat{%
          \vrule\@width\z@\@height\ht\z@
          \vrule\@height\z@\@width\wd\z@}$}%
    \dp\tw@-\ht\z@
    \@tempdima\ht\z@ \advance\@tempdima2\ht\tw@ \divide\@tempdima\thr@@
    \setbox\tw@\hbox{%
       \raise\@tempdima\hbox{\scalebox{1}[-1]{\lower\@tempdima\box
\tw@}}}%
    {\ooalign{\box\tw@ \cr \box\z@}}}
\makeatother

\usepackage{mathtools}
\mathtoolsset{showonlyrefs,showmanualtags}

\usepackage{hyperref} 
\hypersetup{
    colorlinks=true,       
    linkcolor=blue,          
    citecolor=magenta,        
    filecolor=magenta,      
    urlcolor=cyan           
}

\usepackage[msc-links]{amsrefs}

\theoremstyle{plain} 
\newtheorem{lemma}[equation]{Lemma}

\newtheorem{theorem}[equation]{Theorem}

\theoremstyle{definition}
\newtheorem{definition}[equation]{Definition}

\theoremstyle{remark}
\newtheorem{remark}[equation]{Remark}

\numberwithin{equation}{section}

\title[restricted testing condition]{restricted testing conditions for the multilinear maximal operator}
\author[R. J. Chen]{Ruijuan Chen}
\address{Guangling College, Yangzhou University, Yangzhou 225002, China}
\email {rjchen\_yzu@126.com}

\author[W. Chen]{Wei Chen}
\address{School of Mathematical Sciences, Yangzhou University, Yangzhou 225002, China}
\email {weichen@yzu.edu.cn}

\author[F. L. Kong]{Fanliang Kong}
\address{School of Public Education, Xinjiang Institute of Engineering, Urumqi 830000, China}
\email {fanliangkong@cumt.edu.cn}

\thanks{The research of W. Chen is supported by the National Natural Science Foundation of China(11971419, 11771379), the Natural Science Foundation of Jiangsu Province(BK20161326), and the Jiangsu Government Scholarship for Overseas Studies(JS-2017-228). The research of R. J. Chen is supported by College Foundation of Guangling College of Yangzhou University(ZKYB180021)
}

\begin{document}
\begin{abstract}
Restricted testing conditions were considered recently. For the maximal operator, Hyt\"{o}nen, Li and Sawyer \cite{180904873} first obtained parental testing condition.
Later, they \cite{181111032} showed that it suffices to restrict testing to doubling cubes. Chen and Lacey \cite{MR4163998} gave a similar restricted testing condition. In our paper, we discuss a version of the latter in the multilinear setting.
\end{abstract}

	\maketitle
\tableofcontents

\section{Introduction} 

Let $\mathbb R^d$ be the $d\hbox{-dimensional}$ real Euclidean
space and $f$ a real valued measurable function, the Hardy-littlewood maximal function is defined by
$$Mf(x)=\sup\limits_{x\in Q}\frac{1}{|Q|}\int_Q|f(y)|dy,$$
 where $Q$ is a cube with its sides parallel to the coordinate
axes and $|Q|$ is the Lebesgue measure of $Q.$

A weight will be a nonnegative locally integrable function.
Let $u,~v$ be two weights. Muckenhoupt \cite{MR0293384} showed that
$$\left\{
  \begin{array}{ll}
    M:L^p(v,\mathbb{R}^d)\rightarrow L^{p,\infty}(u,\mathbb{R}^d)&\hbox{ iff }(u,v)\in A_p, \hbox{ where } p\geq1;\\
    M:L^p(v,\mathbb{R}^d)\rightarrow L^{p}(v,\mathbb{R}^d)&\hbox{ iff }v\in A_p, \hbox{ where } p>1
  \end{array}
\right.$$
The $A_p$ condition
is geometric, meaning to only involve the weights and not the operators.
Later, Sawyer \cite{MR676801} introduced the test condition $S_p$ and characterized the two weight estimates for the Hardy-Littlewood maximal operator. The classical two weight inequality due to Sawyer \cite{MR676801} is below.
\begin{theorem}\label{t:sawyer} For two weights $ (\omega, \sigma )$ we have the inequality
\begin{equation*}
\lVert M (\sigma f)\rVert _{L ^{p} (\omega)} \lesssim \lVert f\rVert _{L ^{p} (\sigma )}
\end{equation*}
if and only if the testing inequality below holds:
\begin{equation*}
\sup _{Q \;:\; \sigma (Q) >0} \sigma (Q) ^{-1/p} \lVert  \mathbf 1_{Q} M (\sigma \mathbf 1_{Q})\rVert _{L ^{p} (\omega)} < \infty .
\end{equation*}
\end{theorem}

The testing condition essentially amounts
to testing the uniform estimates on characteristic functions of all cubes.
Recent papers Hyt\"{o}nen, Li and Sawyer \cites{180904873,181111032} began a study of a weaker class of
testing inequalities in the two weight setting. Their papers include interesting motivation and background.
They introduced four such conditions in \cites{180904873} and restricted testing to doubling cubes in the two weight norm inequality for $M$ in \cite{181111032}.
Chen and Lacey \cite{MR4163998} gave a similar condition: Test the maximal function on indicators of
cubes $Q$ which have \emph{some} parent on which $ \sigma $ is doubling. We recall the condition and result in \cite{MR4163998}.

\begin{definition}A weight  $\omega$ is a non-negative Borel measure on $ \mathbb R ^d$,
and given two weights $ \omega ,\sigma $  we say that  $ (\omega, \sigma ) \in A _{p}$ if the constant
\begin{equation} \label{e:Ap}
[\omega, \sigma ] _{p} = \sup _{Q}   \langle \omega  \rangle _Q ^{1/p} \langle \sigma  \rangle ^{1/p'}_Q , \qquad  p' = \tfrac p{p-1}.
\end{equation}
where here and throughout $  \langle \omega \rangle_Q = \lvert  Q\rvert ^{-1} \int _{Q} \omega \; dx $.
\end{definition}

\begin{definition}\label{d:parent} Given two weights $ (\omega, \sigma )$, and $ 1< p , \rho , D< \infty $ we say that
$ (\omega, \sigma )$ satisfy a $ (p , \rho , D)$ parent doubling testing condition if there is a positive finite constant $ \mathfrak P = \mathfrak P _{p , \rho, D} = \mathfrak P (\omega, \sigma ,d,  p, \rho ,D)$ so that we have
\begin{equation}\label{e:parent1}
\lVert  \mathbf 1_{Q} M (\sigma \mathbf 1_{Q})\rVert _{L ^{p} (w)}\leq \mathfrak P \sigma (Q) ^{1/p},
\end{equation}
for every cube $ Q$ for which there is a second cube $ P \supset    Q$, with $   \ell (P)  \geq    \rho   \ell (Q)$,
and $ \sigma (P) \leq D \sigma (Q),$ where $ \ell (Q) = \lvert  Q\rvert ^{1/d} $ is the side length of $ Q$.
\end{definition}

\begin{theorem}\label{t:parent} Let   $ 1< p , \rho <\infty$.  There is a constant $ D = D _{d, p , \rho }$ so that
for any pair of weights $ (\omega,\sigma)$ we have
\begin{equation} \label{chen-lacey}
\lVert M (\sigma \cdot )\rVert _{L ^{p} (\sigma ) \to L ^{p} (\omega)}
\simeq    [\omega, \sigma ]_p + \mathfrak P _{ p ,\rho, D}.
\end{equation}
\end{theorem}

The proof of Theorem \ref{t:parent} in \cite{MR4163998} relies on the following essential
ingredients:
\begin{itemize}
  \item classical Sawyer's Theorem \ref{t:sawyer};
  \item $\rho-\hbox{adic}$ grid (dyadic grids in details);
  \item splitting of the subsets of any set in $\rho-\hbox{adic}$ grid (The Testing, The Top and The Small).
\end{itemize}

In our paper, we will give a multilinear version of the Theorem \ref{t:parent}. As far as the authors know, there is no perfect multilinear Sawyer's Theorem. Thus we should find a new ingredient in place of the first one in the above statement.

We state the notation that we will follow in the sequel related to some constants involved in the multiple theory of weights. To define these constants, let $\sigma_1,\ldots,\sigma_m$ and $\omega$ be weights and we denote $\overrightarrow{\sigma}=(\sigma_1,\ldots,\sigma_m)$.  Also let $1<p_1,\ldots,p_m<\infty$ and $p$ be numbers such that $\frac{1}{p}=\frac{1}{p_1}+\dots+\frac{1}{p_m}$ and denote $\overrightarrow p = (p_1,\ldots, p_m)$.

The new multilinear maximal function
\begin{equation*}\label{multi_maximal_operator}\mathcal{M}(\overrightarrow{f})(x)=\mathcal{M}(f_1,...,f_m)(x)= \sup\limits_{x\in Q}
\prod\limits_{i=1}\limits^{m}\frac{1}{|Q|}\int_Q|f_i(y_i)|dy_i, \quad x\in \mathbb R^d
\end{equation*}
associated with cubes with sides parallel to the coordinate
axes was first defined and the corresponding weight theory was studied in \cite{MR2483720}.
The relevant class of multiple weights for $\mathcal{M}$ is given by the condition $A_{\overrightarrow{p}}$
\cite[Definition 3.5]{MR2483720}.

 We recall that $(\omega,\overrightarrow \sigma)$ satisfies the $A_{\overrightarrow P}$ condition if
  \begin{equation}\label{Ap_constant}
  [\omega,\overrightarrow \sigma]_{A_{\overrightarrow P}}:=\sup_{Q}\Big(\frac{1}{|Q|}\int_Q \omega \Big)\prod_{i=1}^m\Big(\frac{1}{|Q|}\int_Q \sigma_i\Big)^{p/p'_i}<\infty.
  \end{equation}

In order to establish the generalization of Sawyer's theorem (Theorem \ref{t:sawyer}) in the multilinear setting, Chen and Dami\'{a}n \cite{MR3118310} introduced
a reverse H\"{o}lder's condition $RH_{\overrightarrow{p}}$ on the weights and established the multilinear version of Sawyer's result with a testing condition $S_{\overrightarrow P}$.
Later on, the condition $RH_{\overrightarrow{p}}$ was used in \cite{MR3424618,MR3350283,MR3750268}. Recently, Cruz-Uribe and Moen \cite{MR4069957} proved a multilinear version of the reverse H\"{o}lder's inequality
in the theory of Muckenhoupt $A_p$ weights.
Note that if $v=\prod_{i=1}^m\omega_i^{{p}/{p_i}},$ then the condition
$(v,\overrightarrow{\omega})\in A_{\overrightarrow{p}}$ implies the reverse H\"{o}lder's condition $RH_{\overrightarrow{p}}$ \cite[Proposition 2.3]{MR3424618}. The reverse H\"{o}lder's condition $RH_{\overrightarrow{p}}$ and the testing condition $S_{\overrightarrow P}$ are below.

We recall that
 $(\omega,\overrightarrow \sigma)$ satisfies the $S_{\overrightarrow{P}}$ condition if
   \begin{equation}\label{Sp_constant}
     [\omega,\overrightarrow \sigma]_{S_{\overrightarrow p}} = \sup_Q \Big(\int_Q \mathcal{M}(\overrightarrow{\sigma\mathbf 1_{Q}})^{p} \omega dx \Big)  \Big(\prod^m_{i=1}\sigma_i(Q)^{\frac{p}{p_i}}\Big)^{-1}<\infty,
   \end{equation}
where $\overrightarrow{\sigma\mathbf 1_{Q}}=(\sigma_1\mathbf 1_{Q},\ldots,\sigma_m\mathbf 1_{Q})$ and all the suprema in the above definitions are taken over all cubes $Q$ in $\mathbb{R}^d$.

We say that $\overrightarrow \sigma$ satisfies the $RH_{\overrightarrow p}$ condition if
there exists a positive constant $C$ such that

  \begin{equation}\label{RH_constant}
     \prod^m_{i=1}\Big(\int_Q\sigma_i dx \Big)^{\frac{p}{p_i}} \leq  C \int_Q\prod^m_{i=1}\sigma_i^{\frac{p}{p_i}} dx.
  \end{equation}
We denote by $[\overrightarrow \sigma]_{RH_{\overrightarrow P}}$ the smallest constant $C$ in \eqref{RH_constant}.

Now we state our restricted testing condition and main result.
\begin{definition}\label{d:parent} Given weights $ (\omega, \overrightarrow{\sigma} )$, and $ 1< \rho , D< \infty $ we say that
$ (\omega, \overrightarrow{\sigma} )$ satisfy a $ (\overrightarrow{p} , \rho , D)$ parent doubling testing condition if there is a positive finite constant $ \overrightarrow{\mathfrak P} = \overrightarrow{\mathfrak P} _{\overrightarrow{p}, \rho , D} = \overrightarrow{\mathfrak P} (\omega, \overrightarrow{\sigma} ,d,  \overrightarrow{p}, \rho ,D)$ so that we have
\begin{equation}\label{e:parent1}
\lVert  \mathbf{1}_Q\mathcal{M}(\overrightarrow{\sigma\mathbf{1}_Q})\rVert^p _{L ^{p} (\omega)}
\leq \overrightarrow{\mathfrak P} \prod^m_{i=1}\sigma_i(Q)^{\frac{p}{p_i}}
\end{equation}
for every cube $ Q$ for which there is a second cube $ P \supset Q$ and an index $1\leq i\leq m$ such that $   \ell P  \geq    \rho   \ell Q$,
and $ \sigma_i (P) \leq D \sigma_i (Q)$.
\end{definition}

\begin{theorem}\label{t:mulparent} Let   $ 1<\rho <\infty$ and let $\overrightarrow \sigma\in RH_{\overrightarrow p}$.  There is a constant $ D = D _{d, \overrightarrow{p} , \rho }$ so that
for weights $ (\omega, \overrightarrow{\sigma})$ we have

\begin{align}
\label{t:ineq1}[\omega,\overrightarrow \sigma]_{A_{\overrightarrow P}} + \overrightarrow{\mathfrak P} _{\overrightarrow{p}, \rho , D}&\lesssim\lVert \mathcal{M} (\overrightarrow{\sigma\cdot}  )\rVert _{\prod^m_{i=1}L ^{p_i} (\sigma_i ) \to L ^{p} (\omega)},\\
\label{t:ineq2}\lVert \mathcal{M} (\overrightarrow{\sigma\cdot}  )\rVert _{\prod^m_{i=1}L ^{p_i} (\sigma_i ) \to L ^{p} (\omega)}
&\lesssim   ([\omega,\overrightarrow \sigma]_{A_{\overrightarrow P}} + \overrightarrow{\mathfrak P} _{\overrightarrow{p}, \rho , D})[\overrightarrow \sigma]_{RH_{\overrightarrow P}},
\end{align}
where $\lVert \mathcal{M} (\overrightarrow{\sigma\cdot}  )\rVert _{\prod^m_{i=1}L ^{p_i} (\sigma_i ) \to L ^{p} (\omega)}$ is the smallest constants such that
\begin{equation}
\label{Th_change}\|\mathcal{M}(\overrightarrow{f\sigma})\|^p_{L^p(\omega)}\leq
C\prod\limits^m_{i=1}\|f_i\|^p_{L^{p_i}(\sigma_i)},
~\forall f_i\in L^{p_i}(\sigma_i).
\end{equation}
\end{theorem}

\begin{remark}\label{remark0}
The proof of \eqref{t:ineq2} is essentially based on four observations:
\begin{itemize}
  \item multilinear version of the Carleson embedding theorem (Lemma \ref{Carleson_lemma_multi});
  \item reverse H\"{o}lder's condition $RH_{\overrightarrow p}$ \eqref{RH_constant}
  \item $\rho-\hbox{adic}$ grid (dyadic grids in details);
  \item splitting of the sparse family in $\rho-\hbox{adic}$ grid (The Testing, The Top and The Small).
\end{itemize}
\end{remark}

\begin{remark}Recall the multilinear Sawyer's Theorem \cite[Theorem 1]{MR3118310}
\begin{equation*}[\omega,\overrightarrow \sigma]_{S_{\overrightarrow p}}\lesssim \lVert \mathcal{M} (\overrightarrow{\sigma\cdot}  )\rVert _{\prod^m_{i=1}L ^{p_i} (\sigma_i ) \to L ^{p} (\omega)}
\lesssim   [\omega,\overrightarrow \sigma]_{S_{\overrightarrow p}}[\overrightarrow \sigma]_{RH_{\overrightarrow P}}.
\end{equation*}
If we use it in place of the Carleson embedding theorem in Remark \ref{remark0}, the exponent of $[\overrightarrow \sigma]_{RH_{\overrightarrow P}}$ is $2.$
\end{remark}

\begin{remark}Let $p_i=q>1,$ $f_i=f$ and $\sigma_i=\sigma$ with $i=1,2,...,m.$ Then
\begin{enumerate}
  \item The reverse condition \eqref{RH_constant} is trivial with $[\overrightarrow \sigma]_{RH_{\overrightarrow P}}=1.$
  \item $[\omega, \sigma ]^q _{q}=[\omega,\overrightarrow \sigma]_{A_{\overrightarrow P}}.$
  \item $\mathfrak P _{q, \rho , D}^q=\overrightarrow{\mathfrak P} _{\overrightarrow{p}, \rho , D}.$
  \item $\lVert M (\sigma \cdot )\rVert^q _{L ^{q} (\sigma ) \to L ^{q} (\omega)}=\lVert \mathcal{M} (\overrightarrow{\sigma\cdot}  )\rVert _{\prod^m_{i=1}L ^{p_i} (\sigma_i ) \to L ^{p} (\omega)}.$
  \item $[\omega,\sigma]^q_{S_{q}}=[\omega,\overrightarrow \sigma]_{S_{\overrightarrow p}}.$
\end{enumerate}

It follows that our main result
\begin{align}
[\omega,\overrightarrow \sigma]_{A_{\overrightarrow P}} + \overrightarrow{\mathfrak P} _{\overrightarrow{p}, \rho , D}&\lesssim\lVert \mathcal{M} (\overrightarrow{\sigma\cdot}  )\rVert _{\prod^m_{i=1}L ^{p_i} (\sigma_i ) \to L ^{p} (\omega)},\\
\lVert \mathcal{M} (\overrightarrow{\sigma\cdot}  )\rVert _{\prod^m_{i=1}L ^{p_i} (\sigma_i ) \to L ^{p} (\omega)}
&\lesssim   ([\omega,\overrightarrow \sigma]_{A_{\overrightarrow P}} + \overrightarrow{\mathfrak P} _{\overrightarrow{p}, \rho , D})[\overrightarrow \sigma]_{RH_{\overrightarrow P}},
\end{align}
reduces to
\begin{equation*}
\lVert M (\sigma \cdot )\rVert _{L ^{p} (\sigma ) \to L ^{p} (\omega)}
\simeq    [\omega, \sigma ]_p + \mathfrak P _{ p ,\rho, D}.
\end{equation*}
which is the main result of \cite{MR4163998}
\end{remark}

\section{Preliminaries}
\label{sec:2}

Before proving our main results, we first recall some definitions and results related to dyadic grids (see \cite[P. 167]{MR3302105} for more information).

Recall that the standard dyadic grid ${\mathcal D}$ in ${\mathbb R}^d$ consists of the cubes

\begin{equation*}
  2^{-k}([0,1)^d+j),\quad k\in{\mathbb Z}, j\in{\mathbb Z}^d.
\end{equation*}

By a general dyadic grid ${\mathcal{D}_{\alpha}}$ we mean a collection of cubes with the following properties:

\begin{enumerate}
  \item for any $Q\in {\mathcal{D}_{\alpha}}$ its sidelength $\ell_Q$ is of the form $2^k, k\in {\mathbb Z}$
  \item $Q\cap R\in\{Q,R,\emptyset\}$ for any $Q,R\in {\mathcal{D}_{\alpha}}$.
  \item the cubes of a fixed sidelength $2^k$ form a partition of ${\mathbb R}^d$.
\end{enumerate}

For $\mathcal{D}_{\alpha},$ we say that $\mathcal{S}_{\alpha}=\{Q_j^k\}\subseteq\mathcal{D}_{\alpha}$ is a {\it sparse family} of cubes if:

\begin{enumerate}
  \item the cubes $Q_j^k$ are disjoint in $j$, with $k$ fixed.
  \item if $\Omega_k=\cup_jQ_j^k$, then $\Omega_{k+1}\subset~\Omega_k$.
  \item $|\Omega_{k+1}\cap Q_j^k|\le \frac{1}{2}|Q_j^k|$.
\end{enumerate}

With each set $Q_j^k\in \mathcal{S}_{\alpha}\cap\Omega_{k},$ we associate the set
\begin{equation} \label{d:sparse}E_j^k=Q_j^k\setminus \Omega_{k+1}.\end{equation}  Observe that the sets $E_j^k$ are pairwise disjoint and $|Q_j^k|\le 2|E_j^k|$. If $\mathcal{D}_{\alpha}=\mathcal{D},$ we simply denote $\mathcal{S}_{\alpha}$ by $\mathcal{S}.$

Next we recall a lemma (\cite[Lemma 3]{MR3118310}). The lemma extends to the multilinear setting a standard formulation of the (dyadic) Carleson embedding theorem proved in \cite{MR3092729} and it will allow us to prove our main results.

\begin{lemma}\cite[Lemma 3]{MR3118310}\label{Carleson_lemma_multi} Let $1 < p_i < \infty$  and $p\in (0,\infty)$ satisfying $\frac{1}{p}=\frac{1}{p_1}+\dots+\frac{1}{p_m}.$ Suppose that the nonnegative numbers $\{a_Q\}_Q$ satisfy

\begin{equation}\label{Carleson_assumption}
 \sum_{Q\subset R} a_Q \leq A \int_R \prod_{i=1}^m \sigma_i^{\frac{p}{p_i}}dx, \, \forall R \in \mathcal{D}
\end{equation}

\noindent where $\sigma_i$ are weights for $i=1,\ldots,m$. Then for all $f_i\in L^{p_i}(\sigma_i)$,

\begin{align*}\label{Carleson}
\sum_{Q\in \mathcal{D}} a_Q \Big(\prod_{i=1}^{m} \frac{1}{\sigma_i(Q)}\int_{Q}f_i(y_i)\sigma_i(y_i)d y_i\Big)^p
\leq A \prod_{i=1}^{m} p'_i \lVert f_i\rVert_{L^{p_i}(\sigma_i)}.
\end{align*}
\end{lemma}

In the sequel we will use the following lemma that could be found in the proof of \cite[Theorem 1.7]{MR3092729}.

\begin{lemma}\label{prhp} There are $2^d$ dyadic grids ${\mathcal{D}}_{i}$ such that for any cube $Q\subset {\mathbb R}^d$ there exists a cube $Q_{i}\in {\mathcal{D}}_{i}$
such that $Q\subset Q_{i}$ and $\ell{Q_{i}}\le 6\ell Q$.
\end{lemma}

\section{Proof of Theorem \ref{t:mulparent}} 

\begin{proof}
It is clear that \eqref{Th_change} implies the $S_{\overrightarrow P}$ condition without using $(\omega,\overrightarrow \sigma)\in RH_{\overrightarrow P}$. Because of $\|\mathcal{M}(\overrightarrow{f\sigma})\|^p_{L^{p,\infty}(\omega)}\leq\|\mathcal{M}(\overrightarrow{f\sigma})\|^p_{L^p(\omega)}$ and \cite[Theorem 3.3]{MR2483720}, \eqref{Th_change} implies the $A_{\overrightarrow P}$ condition.
Then we obtain \eqref{t:ineq1} without the assumption of $\overrightarrow \sigma\in RH_{\overrightarrow P}.$
We now turn to the proof of \eqref{t:ineq2}, which is the main content of the Theorem.

Our theorem only claims that there is a sufficiently large doubling parameter $D$ which
can be used for weights $ (\omega, \overrightarrow{\sigma })$.
 Below, we will consider values of $ 1< \rho \leq 2 $.
 For integers $ \nu = 3, 4 ,\dotsc$, and choices of $ \nu-1< \rho \leq \nu$, the argument proceeds by replacing
 the dyadic grids introduced in Section \ref{sec:2} by $\nu$-ary grids. In fact, modifying the Lemma \ref{prhp}, we may use
$\rho$-ary grids instead of dyadic grids. We omit the details.

By Lemma \ref{prhp}, it suffices to prove \eqref{t:ineq2} for the dyadic maximal operators ${\mathcal M}_{{\mathcal{D}}_{i}}$, where
\begin{equation*}{\mathcal M}_{{\mathcal{D}}_{i}}(\overrightarrow{f})(x)= \sup\limits_{x\in Q\in {\mathcal{D}}_{i}}
\prod\limits_{i=1}\limits^{m}\frac{1}{|Q|}\int_Q|f_i(y_i)|dy_i, \quad x\in \mathbb R^d.
\end{equation*}
Since the proof is independent of the particular dyadic grid, without loss of generality we consider ${\mathcal M}_{\mathcal D}$ taken with respect to the standard dyadic grid ${\mathcal D}.$

Set $D=2^{\frac{2mpd}{mp-1}}$ (see Remark \ref{remark2}). It suffices to  show that under the two weight $ A_{\overrightarrow{p}}$ and  $ (\overrightarrow{p}, 2 , D)$
parent testing condition, the maximal function $ \mathcal{M} _{\mathcal D} (\overrightarrow{\sigma \cdot} )$ is bounded from $ \prod_{i=1}^{m}L ^{p_i} (\sigma_i )$ to $ L ^{p} (\omega)$.

Next we proceed as in the proof of \cite[Theorem 1]{MR3118310}. Then there is a sparse family $\mathcal{S}\subseteq\mathcal{D}$ such that
\begin{align*}
    \int_{\mathbb{R}^d} \mathcal{M}_{\mathcal{D}}(\overrightarrow{f\sigma})^p \omega dx
    \lesssim \sum_{Q\in\mathcal{D}} a_Q \left(\prod_{i=1}^m \frac{1}{\sigma_i(Q)} \int_{Q} |f_i| \sigma_i dy_i \right)^p,
\end{align*}
where $a_Q = \omega(E_Q) \left(\prod_{i=1}^m \frac{\sigma_i(Q)}{|Q|}\right)^p$, if $Q=Q_j^k$ for some $(k,j)$ where $E_Q$ denotes the corresponding set $E_j^k$ associated to $Q_j^k$, and $a_Q=0$ otherwise.

If we apply the Carleson embedding Lemma \ref{Carleson_lemma_multi} to these $a_Q$, we will find the desired result provided that

  \begin{equation}\label{t:test}
    \sum_{Q\subset R} a_Q \lesssim ([\omega, \overrightarrow{\sigma} ] _{A_{\overrightarrow{p}}}  + \overrightarrow{\mathfrak P} _{\overrightarrow{p}, \rho , D})[\overrightarrow \sigma]_{RH_{\overrightarrow P}} \int_R \prod_{i=1}^m \sigma_i^{\frac{p}{p_i}}dx, \, R\in\mathcal{D}.
  \end{equation}

  For $R\in\mathcal{D}$, we denote $\mathcal{S}_R=\{Q\subset R: Q\in \mathcal{S}\}.$ Then we obtain

\begin{align}\label{t:test1}
    \sum_{Q\subset R} a_Q &= \sum_{Q\in\mathcal{S}_R} \omega(E_Q) \left(\prod_{i=1}^m \frac{\sigma_i(Q)}{|Q|}\right)^p.
\end{align}

Partition $ \mathcal{S}_R $ into four subcollections using these definitions.
\begin{itemize}\setlength\itemsep{1em}
\item  (Testing Collection) Let $ \mathcal T ^{\ast} $ be the maximal elements   $ Q\in \mathcal D$
with $ Q\subset R$  so that the testing inequality \eqref{e:parent1} holds.
Set $ \mathcal T_Q = \{P\in \mathcal S \;:\; P\subset Q\}$, for $ Q\in \mathcal T ^{\ast} $.  And set $ \mathcal T = \bigcup _{Q\in \mathcal T ^{\ast} } \mathcal T_Q$.

\item   (The Top)  Let $ \mathcal U = \{Q\in \mathcal{S}_R \setminus \mathcal T \;:\;
2 ^{k}\ell Q \geq \ell R\}$.
We choose $ k$ large enough that $ 2 ^{dmkp} k ^{-2} >1$.
These are the cubes which are close to the top cube  $ R$.

\item  (Small $A_{\overrightarrow{p}}$ Cubes) Let $ \mathcal A $ be those cubes $ Q\in \mathcal{S}_R \setminus (\mathcal T \cup \mathcal U ) $ such that
\begin{equation}\label{e:small}
\Big(\frac{1}{|Q|}\int_Q \omega \Big)\prod_{i=1}^m\Big(\frac{1}{|Q|}\int_Q \sigma_i\Big)^{p/p'_i} \leq \frac{[\omega, \overrightarrow{\sigma} ] _{A_{\overrightarrow{p}}}} {  \phi(\ell   R/  \ell   Q)    } ,
\end{equation}
where $\phi(x)=(\log_2 x)^2$ (see Remark \ref{remark1}).
It is clear the local $A_{\overrightarrow{p}}$ constant at $ Q$ is very small.

\item (Remaining Cubes)  Let $ \mathcal L = \mathcal{S}_R \setminus (\mathcal T \cup \mathcal U \cup \mathcal A)$.
\end{itemize}

We show that  the sum in \eqref{t:test1} over each collection satisfies the testing inequality \eqref{t:test}.
The Testing Collection is very easy:
\begin{align*}
\sum_{Q\in\mathcal{T}} \omega(E_Q) \left(\prod_{i=1}^m \frac{\sigma_i(Q)}{|Q|}\right)^p
& \leq \sum_{Q \in \mathcal T }  \int _{E_Q}  \mathcal{M} (\sigma \mathbf 1_{Q}) ^{p}\omega dx\\
& \leq \overrightarrow{\mathfrak P} _{\overrightarrow{p}, \rho , D} [\overrightarrow \sigma]_{RH_{\overrightarrow P}}
\int_R\prod^m_{i=1}\sigma_i^{\frac{p}{p_i}} dx.
\end{align*}

The Top Collection $ \mathcal U$   has at most  $ 2 ^{1+d(k+1)}$ elements, and we just use the $ A_{\overrightarrow{p}}$ condition
to see that
\begin{align*}
\sum_{Q\in\mathcal{U}} \omega(E_Q) \left(\prod_{i=1}^m \frac{\sigma_i(Q)}{|Q|}\right)^p
&\leq\sum_{Q\in\mathcal{U}} \omega(Q) \left(\prod_{i=1}^m \frac{\sigma_i(Q)}{|Q|}\right)^p\\
& \lesssim  [\omega, \overrightarrow{\sigma} ] _{A_{\overrightarrow{p}}}   [\overrightarrow \sigma]_{RH_{\overrightarrow P}}
\int_R\prod^m_{i=1}\sigma_i^{\frac{p}{p_i}} dx.
\end{align*}
The implied constant depends upon $k$, but that is a fixed integer.

The Small $ A_{\overrightarrow{p}}$ Cubes are also trivially summed up, using the condition in \eqref{e:small}.
\begin{align}
\sum_{Q\in\mathcal{A}} \omega(E_Q) \left(\prod_{i=1}^m \frac{\sigma_i(Q)}{|Q|}\right)^p
&\leq\sum_{Q\in\mathcal{A}} \omega(Q) \left(\prod_{i=1}^m \frac{\sigma_i(Q)}{|Q|}\right)^p\\
& \leq  [\omega, \overrightarrow{\sigma} ] _{A_{\overrightarrow{p}}}   [\overrightarrow \sigma]_{RH_{\overrightarrow P}} \sum _{s> k}\sum _{\substack{\ell R=2^s\ell Q\\Q\in \mathcal A}}  \frac{\int_Q\prod^m_{i=1}\sigma_i^{\frac{p}{p_i}} dx}
 { \phi(  \ell   R/  \ell   Q )  } \\
& =  [\omega, \overrightarrow{\sigma} ] _{A_{\overrightarrow{p}}}  [\overrightarrow \sigma]_{RH_{\overrightarrow P}} \sum _{s> k}\frac{1}
 { \phi(  2^s )  }\left(\sum _{\substack{\ell R=2^s\ell Q\\Q\in \mathcal A}} {\int_Q\prod^m_{i=1}\sigma_i^{\frac{p}{p_i}} dx}
\right) \\
\label{S:con}&\lesssim [\omega, \overrightarrow{\sigma} ] _{A_{\overrightarrow{p}}}   [\overrightarrow \sigma]_{RH_{\overrightarrow P}}
\int_R\prod^m_{i=1}\sigma_i^{\frac{p}{p_i}} dx.
\end{align}

Thus, the core of the argument is control of the Remaining Cubes, $ \mathcal L$. We claim that this
collection is empty.

Suppose $ \mathcal L \neq \emptyset $. Thus, there is a cube $ Q \subset R$,
which satisfies $ \ell Q < 2 ^{-k} \ell R $,   fails \eqref{e:small}, and \emph{no ancestor of $ Q$ also  contained inside of $ R$, has a doubling parent.}  The last condition is  very strong.

Let $ Q ^{(1)}$ be the $ \mathcal D$-parent of $ Q$, and let $ Q ^{(j+1)} = (Q ^{(j)}) ^{(1)}$.
Define integer $ n$ by  $ R = Q ^{(n)}$.
For any integer $ 0\leq j < n$ and $1\leq i\leq m$, we
necessarily have $   \sigma_i (Q ^{(j+1)})> D  \sigma_i (Q^{(j)})$, since $ Q ^{(j+1)}$ is a $ \rho $-parent of
$ Q ^{(j)}$.
That is, $ \sigma_i (R) \geq D ^{n} \sigma_i (Q),$ $1\leq i\leq m$.
From this, we see that $n$ cannot be very large.
\begin{align}
[\omega, \overrightarrow{\sigma} ] _{A_{\overrightarrow{p}}}   & \geq \Big(\frac{1}{|R|}\int_R \omega \Big)\prod_{i=1}^m\Big(\frac{1}{|R|}\int_R \sigma_i\Big)^{p/p'_i}\\
& \geq   D ^{n(mp-1)} \Big(\frac{1}{|Q ^{(n)}|}\int_{Q} \omega \Big)\prod_{i=1}^m\Big(\frac{1}{|Q ^{(n)}|}\int_{Q} \sigma_i\Big)^{p/p'_i}
\\
& \geq 2^{-dn}\Big(\frac{D}{2^d}\Big) ^{n(mp-1)} \Big(\frac{1}{|Q |}\int_{Q} \omega \Big)\prod_{i=1}^m\Big(\frac{1}{|Q|}\int_{Q} \sigma_i\Big)^{p/p'_i}\\
& \geq 2^{-dn}\Big(\frac{D}{2^d}\Big) ^{n(mp-1)}\frac{[\omega, \overrightarrow{\sigma} ] _{A_{\overrightarrow{p}}}} {  \phi(  \ell   R/  \ell   Q )   }\\
\label{d:pow}& = [\omega, \overrightarrow{\sigma} ] _{A_{\overrightarrow{p}}}  2^{-dn}\Big(\frac{D}{2^d}\Big) ^{n(mp-1)}n^{-2}\\
\label{d:fixD}& = [\omega, \overrightarrow{\sigma} ] _{A_{\overrightarrow{p}}}  2^{dmnp}n^{-2}.
\end{align}
Note that we have used $D=2^{\frac{2mpd}{mp-1}}$ in \eqref{d:fixD}.
Recall that we choose $k$ large enough that $2 ^{dmkp} k ^{-2} >1.$
We see that $ n < k$.  That is, the cube is in the collection $ \mathcal U$, which is a contradiction.
\end{proof}

\begin{remark}\label{remark1}
Let $q>1.$ We can define $\phi(x)=(\log_2 x)^q$ in our proof. Then our proof is still valid. In fact, we have
$\sum _{s> k}\frac{1} { \phi(  2^s )  }=\sum _{s> k}\frac{1}{s^q}<\infty$ in \eqref{S:con} and $n^{-q}$ in place of $n^{-2}$ in \eqref{d:pow}.
We mention that we do not try to find the optimal $\phi.$
\end{remark}

\begin{remark}\label{remark2}
To determine $D,$ we can let $D=2^{dt}.$ Then we have $2^{dn\left((t-1)(mp-1)-1\right)}$ in place of $2^{dmnp}$ in \eqref{d:fixD}.
It suffices to choose $t$ such that $(t-1)(mp-1)-1>0.$ We do not try to find the optimal $D.$
\end{remark}

\section{Acknowledgement}
We thank the referees for many valuable comments and suggestions. These greatly improve the presentation of our results.

\bibliographystyle{alpha,amsplain}	


\end{document}